\documentclass[conference]{IEEEtran}
\IEEEoverridecommandlockouts

\usepackage{amsfonts,amsmath,graphicx,subfigure}
\usepackage{epstopdf}
\usepackage{siunitx}
\usepackage[hyphens]{url}

\usepackage{cite}
\usepackage{amssymb}
\usepackage{algorithmic}
\usepackage{graphicx}
\usepackage{textcomp}

\newtheorem{example}{Example}[section]

\graphicspath{{./plots/}}

\begin{document}
\Urlmuskip=0mu plus 2mu

\title{Polynomial Preconditioned GMRES to Reduce Communication in Parallel Computing\\
\thanks{Sandia National Laboratories is a multimission laboratory managed and operated by National Technology and Engineering Solutions of Sandia, LLC., a wholly owned subsidiary of Honeywell International, Inc., for the U.S. Department of Energy’s National Nuclear Security Administration under contract DE-NA-0003525.}
}

\author{\IEEEauthorblockN{Jennifer A. Loe}
\IEEEauthorblockA{\textit{Department of Mathematics} \\
\textit{Baylor University}\\
Waco, Texas, USA \\
jennifer\_loe@baylor.edu}
\and
\IEEEauthorblockN{Heidi K. Thornquist}
\IEEEauthorblockA{\textit{Electrical Models and Simulation} \\
\textit{Sandia National Laboratories}\\
Albuquerque, New Mexico, USA \\
hkthorn@sandia.gov}
\and
\IEEEauthorblockN{Erik G. Boman}
\IEEEauthorblockA{\textit{Center for Computing Research} \\
\textit{Sandia National Laboratories}\\
Albuquerque, New Mexico, USA\\
egboman@sandia.gov}
}

\maketitle

\begin{abstract}

 Polynomial preconditioning with the GMRES minimal residual polynomial has the potential to greatly reduce orthogonalization costs, making it useful for communication reduction.  We implement polynomial preconditioning in the Belos package from Trilinos and show how it can be effective in both serial and parallel implementations.   
We further show it is a communication-avoiding technique and is a viable option to CA-GMRES for large-scale parallel computing.
\end{abstract}


\section{Introduction}
As computational models become more complex, efficient analysis of such models will require advanced numerical algorithms that target high-performance computing.  At the core of these analyses is often the solution of large, sparse linear systems $Ax=b$.  The most common algorithms for performing such solves are Krylov subspace methods, like the Generalized Minimum Residual Method (GMRES) \cite{JAL:SaadSchultz}.  GMRES is still considered one of the best general purpose, iterative methods for non-Hermitian linear systems, even though it requires full orthogonalization of the Krylov subspace, $\mathcal{K}_m(A,b)= span\{b, Ab, A^2b, \ldots, A^{m-1}b\}$.  To accelerate the convergence of GMRES, and reduce costs like orthogonalization, it is common to solve a preconditioned linear system.  The preconditioner, $M$, is ideally an inexpensive approximation of $A^{-1}$ and can be applied to the left, $MAx=MB$, or right, $AMy=b$ where $x=My$, of the original linear system.  

Modern architectures often gain computing power through the addition of more processors (cores), rather than through increases in processor clock speed.  To efficiently utilize these architectures, it is necessary to leverage parallel algorithms.  Due to the high communication and global synchronization requirements of Krylov subspace methods, like GMRES, much research has been devoted to communication-avoiding algorithms \cite{JAL:HoemmenThesis, JAL:MinComm, JAL:AvoidComm}.  Among these algorithms are $k$-step methods, which compute several new Krylov subspace vectors between each communication-intensive orthogonalization. 

Polynomial preconditioning is another approach that has potential to reduce communication costs in Krylov solvers. Instead of performing many iterations at once, like the $k$-step methods, this approach packs more work into each iteration.  With polynomial preconditioning, more matrix-vector products are used to form each basis vector before orthogonalizing, giving GMRES more power and potential for convergence for roughly the same amount of communication.   Chebyshev and least-squares polynomials are common choices for polynomial preconditioners since they can be created so that their norms are minimized over a given interval.  Unfortunately, obtaining an effective Chebyshev or least squares polynomial requires estimates for the extreme eigenvalues of the associated matrix, which are often hard to obtain in practice. This is especially true for non-Hermitian linear systems where determining the polynomial may require an estimate of the convex hull of the spectrum \cite[p. 403]{JAL:SaadItMeth}.

Trilinos \cite{JAL:Trilinos} is one of the major software libraries that offers a collection of advanced numerical algorithms for parallel computing, including sparse linear solvers and preconditioners.  While polynomial preconditioners have long been considered for parallel computing (\hspace{1sp}\cite{SaadPractical,LiangThesis} and references therein), they are rarely included in high-performance numerical software collections.  In Trilinos, for example, the only polynomial preconditioners are the Chebyshev and least-squares preconditioners in the IFPACK package and they only work for Hermitian matrices.  The more recent package, IFPACK2, only provides Chebyshev as a smoother for multigrid.  This lack of polynomial preconditioners in high-performance software may be due to the computational cost required to compute eigenvalues for obtaining an effective polynomial.

The GMRES minimum residual polynomial also has small norm over the spectrum of the given matrix, making it an effective preconditioner \cite{JAL:PPGMRES}.  Its construction does not require any explicit information about the spectrum, making it cheaper to obtain than other polynomials.   This preconditioner has been used effectively in several contexts: In \cite{JAL:Thorn, PPArnoldi}, it is used to spectrally transform an operator to solve for eigenvalues.  The work \cite{JAL:PPGMRES} uses the GMRES minimum residual polynomial to precondition serial implementations of GMRES and GMRES-DR, and \cite{JAL:PPBICG} studies the preconditioner for methods IDR(s) and BiCGStab. 


In this paper, we demonstrate an implementation of the GMRES minimum residual polynomial preconditioner using Trilinos.  We discuss issues that arise when implementing the preconditioner into existing software, and we test its potential to accelerate the solution of large-scale linear systems. We illustrate that this preconditioner can be effective in a high-performance computing environment for the following reasons:  

\begin{enumerate}
\item \textbf{Simple and Effective:} The GMRES minimum residual polynomial preconditioner is simple to implement, as discussed in Section \ref{sec:implementation}.  Since no factorization of the matrix is required, this preconditioner is suitable for problems where the full matrix is unavailable.  The effectiveness of the preconditioner is demonstrated in Sections \ref{sec:serialresults}, \ref{sec:parallelresults}.
\item \textbf{Avoiding Communication:}  Applying the polynomial preconditioner requires only sparse matrix-vector products (SpMVs) and vector updates.  This greatly reduces global communication and synchronization from inner products, which may become a bottleneck.  We discuss in Section \ref{sec:CA-Alg} potential to combine the preconditioner with other communication-avoiding kernels. 
\item \textbf{Accelerates Existing Preconditioners:}  The polynomial preconditioner can be combined with existing preconditioners.  We use it to accelerate ILU in Section \ref{sec:parallelresults}.
\item \textbf{Potential for Automation:}  We give a heuristic for automating the selection of the polynomial degree in Section \ref{sec:degreeselection}.  
\end{enumerate}
Because of these advantages, we propose that polynomial preconditioning should be considered as an addition to high-performance solvers software libraries.

\section{Implementing the Minimum Residual Polynomial Preconditioner}
\label{sec:implementation}
Given a linear system $Ax=b$, we obtain the polynomial preconditioner $p(A)$ of degree $deg$. We first build a power basis $V = [v_0, Av_0, \ldots, A^{deg}v_0]$, where $v_0$ is an arbitrary vector.  Then we solve the normal equations 
\begin{equation} 
(AV)^*AVy=(AV)^*v_0. 
\label{powereqn} 
\end{equation}    
The elements of $y$ are the coefficients of $p(A)$, that is, 
\begin{equation}
\label{eq:poly}
 p(A) = y_{deg+1} A^{deg} + y_{deg} A^{deg-1} + \cdots + y_2 A + y_1 . 
\end{equation}
This method becomes unstable as the the columns of $V$ lose linear independence, but its results are generally sufficient for low-degree polynomials.  

Notice that the coefficients of the polynomial depend on the choice of $v_0$.  Ref. \cite{JAL:PPGMRES} suggests using $v_0=b$, the problem right-hand side.  We discuss in Section \ref{sec:choosingvector} why a random vector may instead be preferable.  Note also that the polynomial preconditioner can be easily combined with other preconditioners:  Given a preconditioned system $MAx=Mb$, we can use the operator $MA$ to form the power basis $V$ and obtain a polynomial preconditioner for the already preconditioned system.   

The spectrum of the preconditioned operator $Ap(A)$ will typically be better for convergence of GMRES than that of the original matrix $A$: the small eigenvalues of $A$ will be mapped to well-separated eigenvalues of $Ap(A)$, and other eigenvalues will be clustered near $1$.  Further details on the derivation and algebraic properties of the polynomial, as well as algorithms for more stable implementation, can be found in \cite{JAL:PPGMRES}.

The results presented in this paper are obtained using an implementation of the minimum residual polynomial preconditioner written directly in Belos, the next-generation iterative linear solvers package of Trilinos \cite{JAL:BelosAmesos2}.  While most preconditioners are found in other packages within Trilinos, the Belos package was the most convenient location since it already has code for generating Krylov subspaces and performing orthogonalization.  Furthermore, Belos is designed so that it can be used to precondition itself as an inner-outer solver.  Finally, this implementation provides the advantage that it works for both Epetra and Tpetra-based linear algebra.
 
The GMRES minimum residual polynomial operator is generated in Belos using the vector (MultiVecTraits) and operator (OperatorTraits) abstractions.  The GmresPolyPrec class takes a Belos linear problem and a polynomial degree as arguments to the constructor and computes the coefficients of the polynomial.  It forms the matrices for (\ref{powereqn}) and then uses the LAPACK routines \cite{JAL:LAPACK}, POTRF/POTRS, to compute a Cholesky factorization and solve the system.  The class also provides a function to apply $p(A)$ to a given vector $x$. 

A new Belos SolverManager class, GmresPolyPrecSolMgr, allows current Belos linear solvers to interface to the polynomial preconditioner as they would another linear solver.  The Belos SolverFactory was extended to include the new GmresPolyPrecSolMgr, allowing us to set the polynomial preconditioner as an inner solver and GMRES as an outer solver.  For the results in this paper, we used the Epetra linear algebra interfaces.

\section{Serial Numerical Results}
\label{sec:serialresults}
This section contains results from a serial build of the Belos polynomial preconditioner.  The outer solver is GMRES$(m)$ with two steps of Classical Gram Schmidt (ICGS) orthogonalization.  Each orthogonalization step requires two block inner products and one norm, which we count as three "dot" products per iteration in the figures to follow.  (This is more effective for avoiding communication than modified Gram-Schmidt orthogonalization.)  The algorithm needs $deg+1$ SpMVs and two block inner products to create the polynomial, which is used as a right preconditioner.  We require a relative residual norm less than \num{1e-8} for convergence.

For all results in this section, we use the matrix e20r0100 from Matrix Market \cite{JAL:MatrixMarket}, a real non-Hermitian matrix of size $n=4241$.  This matrix has high condition number, estimated by Matrix Market at \num{2.15e10}, and proves to be extremely difficult for GMRES. In practice, e20r0100 may be best addressed using a direct solver, but it is representative of the difficulties that one may encounter in GMRES.

\begin{example}
\label{ex:serial1}
\normalfont
In this example, we choose subspace size $m=50$ and a right-hand side, $b$, with random entries. We generate the polynomial using $v_0 = b$.  Fig. \ref{JAL:fig2} shows residual norm convergence for unpreconditioned GMRES (indicated by $deg=0$) and polynomial preconditioned GMRES for $deg =3, 5, 7, 10$.  Convergence in relation to SpMVs is shown on the top and convergence relative to inner products is shown on the bottom.   

\begin{figure}
\centering
\includegraphics[scale=.6]{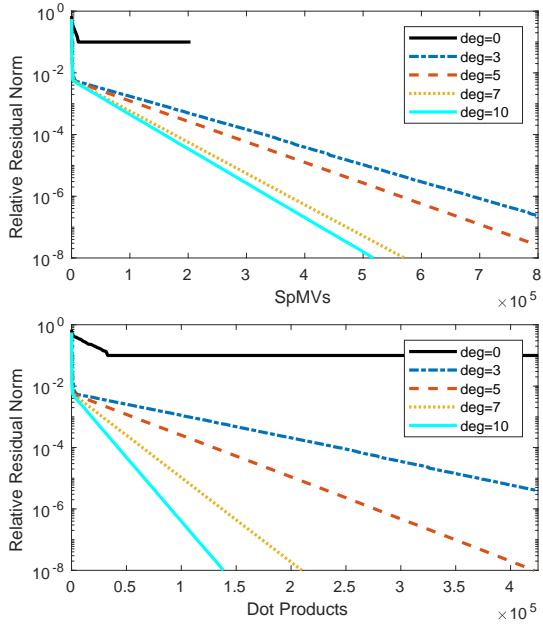}
\caption{Residual norm convergence for the matrix e20r0100 with a random right-hand side.  Subspace size = 50.  Degree $0$ indicates no preconditioning.  All tests were run to $200000$ max iterations.}
\label{JAL:fig2}
\end{figure}
The results, illustrated in Figure \ref{JAL:fig2}, show that without preconditioning the relative residual norm only improves by one order of magnitude before stalling out.  For these experiments we see a distinct improvement when raising the degree of the polynomial.   While the degree $5$ problem converges in $859993$ SpMVs, the polynomial of degree $10$ gives the most improvement, converging in $517452$ SpMVs.  The decrease in inner (dot) products is more substantial.  The degree $5$ problem converges in $421559$ dot products, while the degree $10$ problem only requires $138350$.  This is over three times fewer dot products, and a stark difference from the unpreconditioned problem which stagnated.  

Considering the spectra of $A$ and $Ap(A)$ helps to explain the improvement from polynomial preconditioning.  Figure \ref{JAL:evals} shows the eigenvalues of $A$ (top) and the new spectrum after applying a degree $7$ preconditioner (bottom).  The matrix $A$ is indefinite, having $1199$ eigenvalues that lie in the left half of the complex plane.  The largest eigenvalue with negative real part has magnitude $0.0013$, and the smallest has magnitude $\num{1.4e-06}$.  While the polynomials do not create a significant change in the eigenvalues that have negative real part, they do help to cluster eigenvalues on the right-half plane away from the origin.  We find that the ratio of smallest to largest magnitudes of eigenvalues on the right half plane is $0.0028$ without preconditioning.  With only a degree $3$ polynomial, this ratio improves almost ten times to $0.0255$.  When the polynomial is degree $10$, the ratio has improved to $0.0807$.  

\begin{figure}
\centering
\subfigure[Eigenvalues of $A$]{\includegraphics[scale=.5]
{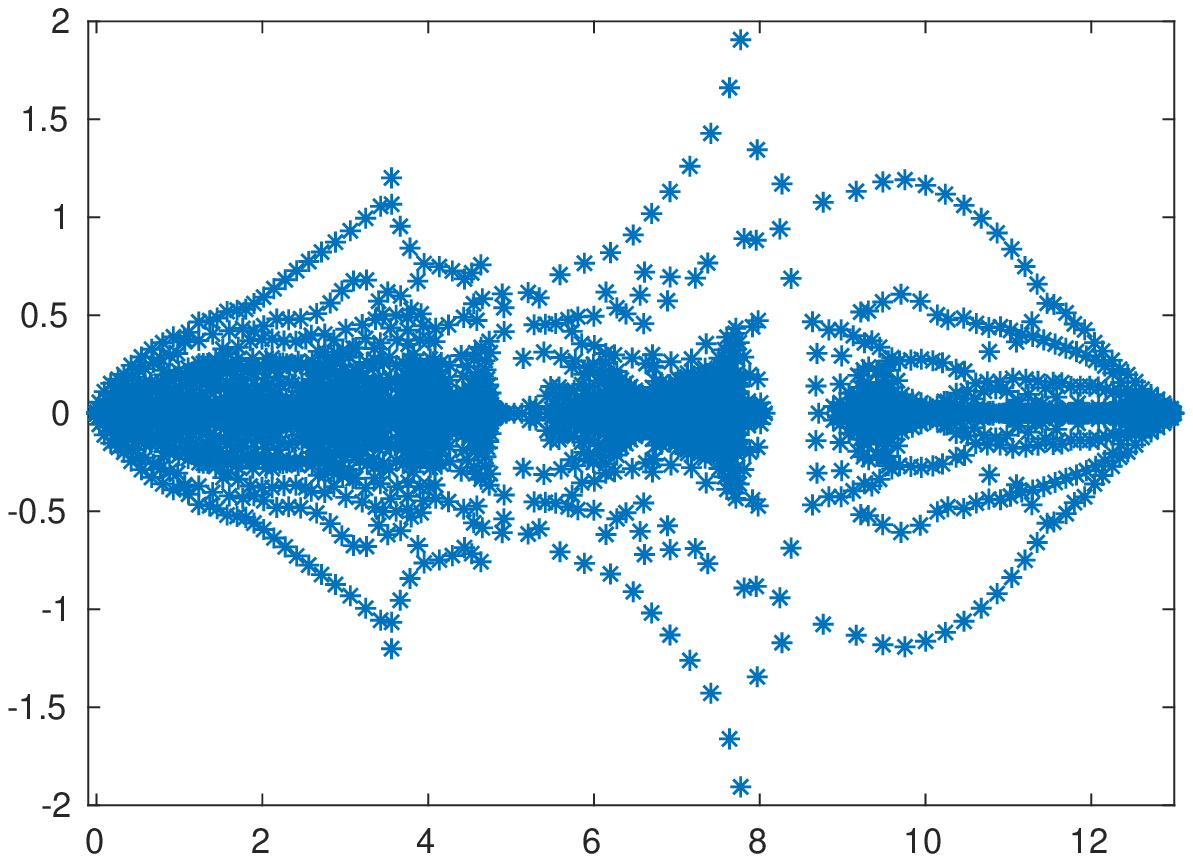}}
\subfigure[Eigenvalues of $Ap(A)$. Degree of $p$ is $7$]{\includegraphics[scale=.5]
{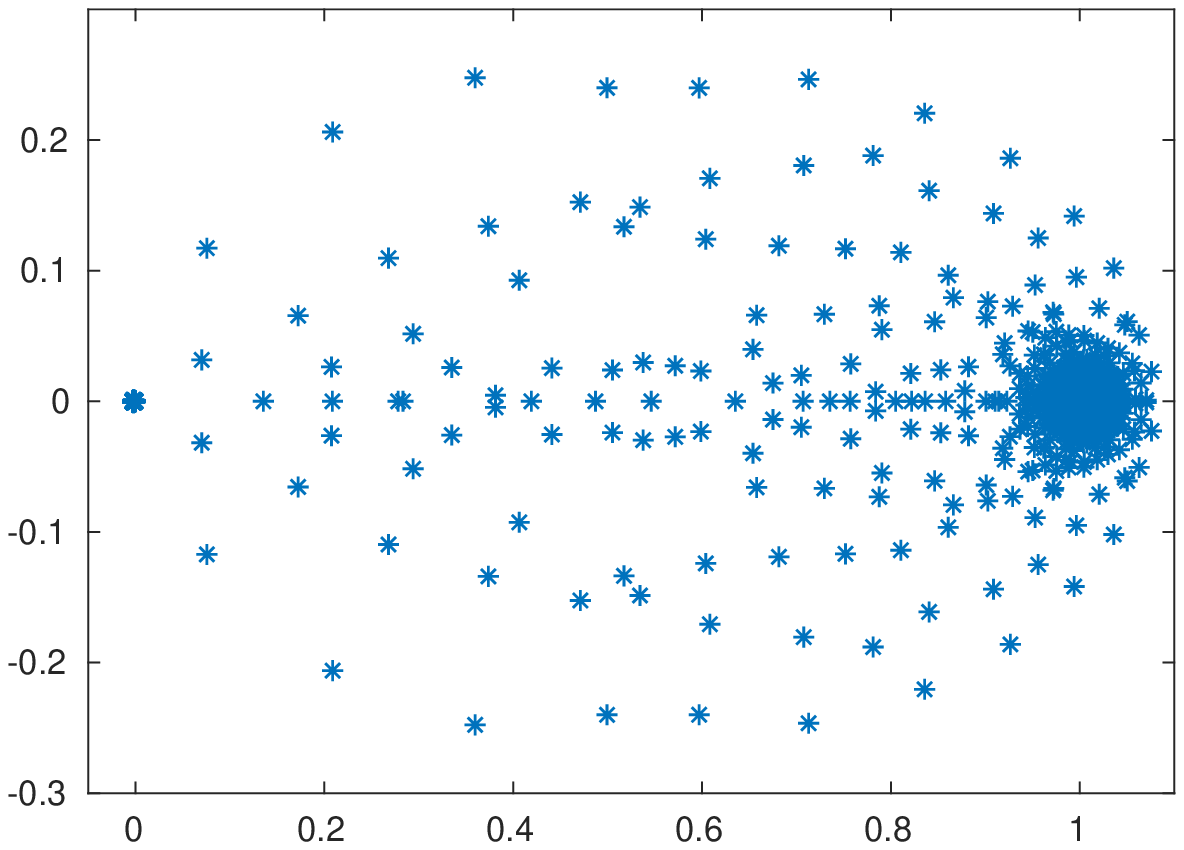}}
\caption{Eigenvalues of the preconditioned matrix $Ap(A)$ where $A$ is the matrix e20r0100.  As the degree of $p(A)$ increases, the eigenvalues are more clustered around $1$ and more loosely scattered near zero.}
\label{JAL:evals}
\end{figure}
\end{example}

\begin{example}
\normalfont
In our next example, we use the same problem and polynomials from Example \ref{ex:serial1}, but we increase the subspace size to $m=100$.  Though GMRES no longer stalls with the larger subspace, convergence is too slow to run to completion.  We estimate that $1335530$ SpMVs and $3966900$ dot products are needed to converge.  Figure \ref{JAL:fig3} shows the improvement with preconditioning.  The degree $3$ preconditioned problem converges faster, but still too slowly to run to completion. The polynomial of degree $5$ helps to attain convergence at a cost of $147421$ SpMVs and $72974$ dot products.  With degree $10$, the cost is $40008$ SpMVs and $10799$ dot products.  Thus, if this computation was run in parallel, we would reduce global communication calls by about $2.5$ orders of magnitude over no preconditioning.  With that comes approximately $1.5$ orders of magnitude improvement in matrix-vector products, reducing processor-to-processor communication as well. 
\begin{figure}
\centering
\includegraphics[scale=.6]{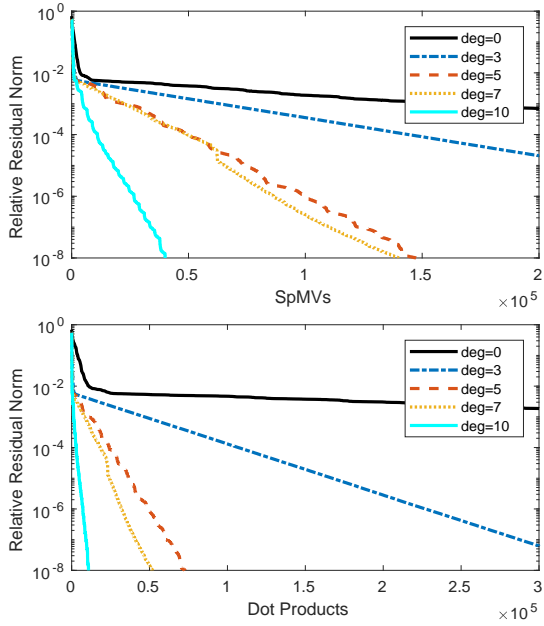}
\caption{Residual norm convergence in terms of matrix-vector products for the matrix e20r0100 with a random right-hand side.  Subspace size $= 100$.  Deg $0$ indicates no preconditioning.}
\label{JAL:fig3}
\end{figure}
 
\end{example}

\section{A Degree Selection Strategy}
\label{sec:degreeselection}
\normalfont
It may be difficult for the user to determine when it is best to stop raising the polynomial degree.  Raising the degree often results in a better preconditioner, but it can reach a point of diminishing returns.  The polynomial preconditioner can decrease both the number of inner products and SpMVs required to converge, but sometimes inner products are reduced at the expense of more SpMVs.  Fortunately, many matrices are stored so that SpMVs only require communication with neighboring processors.  
Thus, for communication reduction, it may be more beneficial to perform extra SpMVs in order to avoid operations that require synchronous global communication, like inner products.

We ran serial tests on several different matrices to determine the effects of raising the polynomial degree.  All tests were performed with a right-hand side $b=Ax$ where $x$ is a randomly generated solution vector.  We let $v_0=b$ and choose a maximum subspace size of $50$.  Matrices bwm2000, orsirr1, s1rmq4m1, and e20r0100 can be obtained via Matrix Market. The matrix BiDiag1, 
with $n=2000$, has $1,2,\ldots, 2000$ on the diagonal and $0.05$ on all elements of the superdiagonal. 
Matrix BiDiag2, with $n=5000$, 
has $0.1, 0.2, \ldots, 0.9, 1, 2, \ldots, 4991$ on the diagonal and $0.2$ on all elements of the superdiagonal.   

\begin{figure}
\centering
\includegraphics[scale=.6]{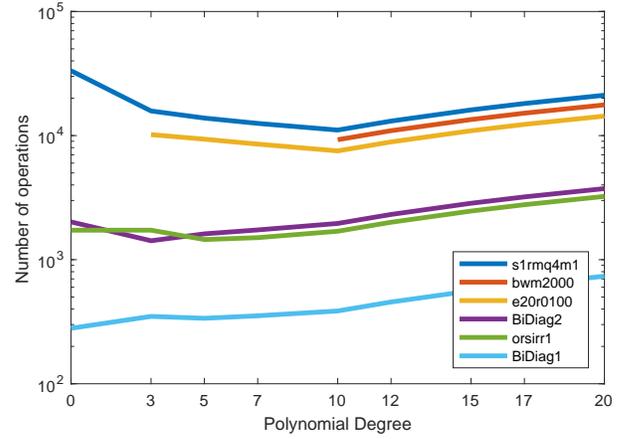}
\caption{The number of SpMVs required to reach convergence for several preconditioned matrices with different polynomial degrees.  Subspace size is $50$.}
\label{fig5}
\end{figure}

\begin{figure}
\centering
\includegraphics[scale=.6]{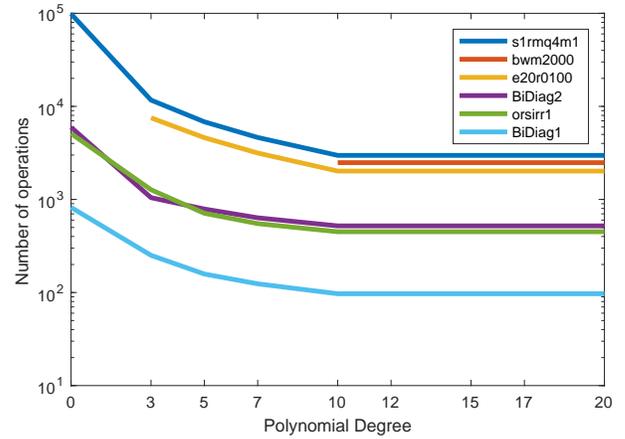}
\caption{Total number of dot products (inner products plus norms) from orthogonalization for several preconditioned matrices with different polynomial degrees.  Subspace size is $50$.}
\label{fig7}
\end{figure}

Figure \ref{fig5} shows the number of SpMVs required to reach a relative residual tolerance of $\num{1e-8}$ for polynomials of degrees $3,5,7,10,12,15,17,20$.  Results for no preconditioning correspond to degree $0$ on the plot.  Figure \ref{fig7} shows the corresponding number of inner products required for convergence, where one inner product is counted for each of two passes of Gram-Schmidt orthogonalization and one more for the norm. 
For matrices bwm2000 and e20r0100, convergence stagnates with no preconditioning.  The problem e20r0100 first converges with the preconditioner of degree $3$, and bwm2000 first converges with degree $10$.

The results suggest that preconditioning is more likely to reduce matrix-vector products for difficult problems than for simpler ones.  For the easiest problem, BiDiag1, matrix-vector products increase with preconditioning, even for degree $3$.  For matrices e20r0100, bwm2000, and s1rmq4m1, which were most difficult, the expense of SpMVs decreases with preconditioning up until degree $10$. For the other two problems, the number of SpMVs decreases slightly for very low-degree polynomials and begins to rise again, with no savings after degree $10$.  

Unlike with SpMVs, polynomial preconditioning is consistent in reducing the number of inner products for all problems in Figure \ref{fig7}, regardless of difficulty.  By the time the polynomial degree is increased to $10$, the number of inner products has decreased by approximately an order of magnitude or more for all problems.  However, after degree $10$, the number of inner products remains constant while the number of SpMVs is increasing.  

We found that after the polynomial degree gets large enough, increasing it failed to result in new coefficients of significant magnitude.   Example coefficients for the matrix s1rmq4m1 are shown in Table \ref{coefftable}.  Notice that for degrees $10, 12$, and $15$, the first eleven polynomial coefficients remain the same.  The additional coefficients in the polynomials of degrees $12$ and $15$ are so near zero that they do not provide any additional information.  This explains why the degree $12$ and degree $15$ preconditioners give no improvement in cost over the degree $10$ preconditioner.  In fact, they are more expensive due to the extra SpMVs incurred with a near-zero coefficient.  Polynomial coefficients for the other matrices tested followed a similar trend, becoming very small at high degrees.  

\begin{table}[]
\centering
\caption{Coefficients for the polynomial $p(A)$ generated with the matrix s1rmq4m1.}
\resizebox{.45\textwidth}{!}{%
\begin{tabular}{|r|r|r|r|}
\hline
\multicolumn{1}{|c|}{\textbf{Deg 7}} & \multicolumn{1}{c|}{\textbf{Deg 10}} & \multicolumn{1}{c|}{\textbf{Deg 12}} & \multicolumn{1}{c|}{\textbf{Deg 15}} \\ \hline
3.70798e-05                          & 6.08343e-05                          & 6.08343e-05                          & 6.08343e-05                          \\ \hline
-5.2613e-10                          & -1.45759e-09                         & -1.45759e-09                         & -1.45759e-09                         \\ \hline
3.82756e-15                          & 1.86807e-14                          & 1.86807e-14                          & 1.86807e-14                          \\ \hline
-1.59154e-20                         & -1.45193e-19                         & -1.45193e-19                         & -1.45193e-19                         \\ \hline
3.93091e-26                          & 7.29453e-25                          & 7.29453e-25                          & 7.29453e-25                          \\ \hline
-5.69776e-32                         & -2.44543e-30                         & -2.44543e-30                         & -2.44543e-30                         \\ \hline
4.47271e-38                          & 5.52107e-36                          & 5.52107e-36                          & 5.52107e-36                          \\ \hline
-1.46677e-44                         & -8.28868e-42                         & -8.28868e-42                         & -8.28868e-42                         \\ \hline
                                     & 7.92937e-48                          & 7.92937e-48                          & 7.92937e-48                          \\ \hline
                                     & -4.3729e-54                          & -4.3729e-54                          & -4.3729e-54                          \\ \hline
                                     & 1.05777e-60                          & 1.05777e-60                          & 1.05777e-60                          \\ \hline
                                     &                                      & 7.47962e-208                         & 7.47962e-208                         \\ \hline
                                     &                                      & 6.051e-237                           & 6.051e-237                           \\ \hline
                                     &                                      &                                      & -3.85124e-257                        \\ \hline
                                     &                                      &                                      & -1.25918e-274                        \\ \hline
\end{tabular}
}
\label{coefftable}
\end{table}

This effect may also be due to the ill-conditioned problem of computing coefficients via the normal equations with a power basis.  We observe that the appearance of near-zero coefficients corresponds with a positive return value \textit{info} in the LAPACK function POTRF when forming the polynomial, which means that the matrix [$(AV)^*(AV)$, in our case] is not positive definite. 
Of the six matrices discussed in this section, three first give positive return values starting with degree $10$, and the other three examples begin to give warnings at degree $12$.  Despite this warning, all coefficients are still computed.  In other examples, NaNs were computed after the LAPACK error occurred and the polynomial degree was raised too high. 

 We also tried using a more stable LAPACK routine POSVX, which equilibrates the system before Cholesky factorization and/or improves the solution using iterative refinement.  Neither of these options resulted in better polynomial coefficients.  Another option is to find a QR factorization for solving the normal equations.  In \cite[p. 11]{JAL:PPGMRES}, this resulted in less accurate polynomial coefficients.  Thus, we do not consider it here. 

We conjecture that the best polynomial preconditioner constructed with the power basis method will have the highest degree possible without a warning from LAPACK.  Based on the examples above, this polynomial seems likely to minimize the number of inner products and norms while avoiding extra SpMVs.  This strategy can be easily implemented for automatic degree selection. 

It is worth noting that there are examples, such as Sherman5 from Matrix Market, where this degree selection strategy fails.  This indefinite matrix is an extremely difficult problem for GMRES.  The polynomial of degree $7$ was a very successful preconditioner because the spectrum of $Ap(A)$ was entirely in the one side of the complex plane.   With higher degree polynomials, the matrix $Ap(A)$ was once again indefinite and GMRES did not converge, but LAPACK did not give positive return values until degree $15$.  In such instances, it may be best to take the auto-selection degree as an upper bound and try to obtain results with lower-degree polynomials.  

\section{Choosing a Vector to Generate the Polynomial}
\label{sec:choosingvector}
All experiments thus far have successfully generated the polynomial preconditioner using $v_0 = b$, the problem right-hand side.  This choice worked well in the previous sections because the problem right-hand side was generated using randomization.  More structured right-hand sides may generate a poor polynomial preconditioner.

\begin{figure}[]
\centering
\subfigure[Polynomial from bad start vector]{\includegraphics[scale=.6]{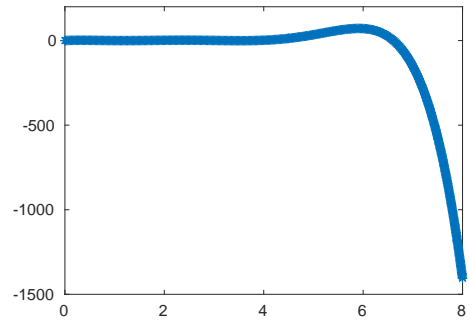}}
\hfill
\subfigure[Closeup of bad polynomial]{\includegraphics[scale=.6]{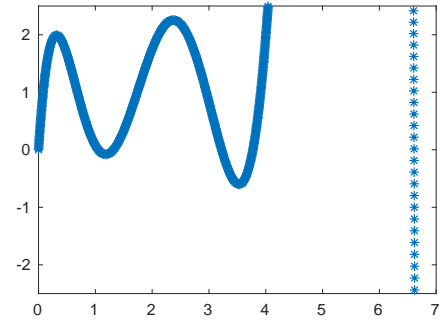}}
\caption{Polynomial $\alpha p(\alpha)$ of degree $5$ for the Laplacian matrix where $v_0 =b$.  Plotted at points between $[0,8]$.  Closeup on the bottom.}
\label{fig:BadPolyLapl}
\end{figure}

Consider the discretized Laplacian equation $-\nabla ^2 u = f$ over a square domain, with constant source function $f(x)\equiv 1$ and zero boundary conditions (Example $1.1.1$ \cite{ElmanFiniteElts}).  The matrix size is $n=40401$.  The eigenvalues of this matrix are all real-valued and lie in the interval $[0,8]$, with several eigenvalues very close to $8$.  All values of the right-hand side vector $b$ are very close to $0$ or $1$.  Figure \ref{fig:BadPolyLapl} shows the polynomial $\alpha p(\alpha)$ of degree $5$ generated with $v_0 = b$.    
The $x$-axis corresponds to the spectrum of $A$, and the $y$-axis shows the range of eigenvalues of $Ap(A)$.  Recall that if $p(A)$ is a good preconditioner, the large eigenvalues of $A$ will be mapped close to $1$ and the small eigenvalues of $A$ will be well-separated between $0$ and $1$.  This polynomial does nothing of the sort.  The largest eigenvalues near $8$ are mapped to near $-1400$, and the eigenvalues in the middle of the spectrum are mapped to values as small as $-1/2$ up to larger than $2$. 

The preconditioned matrix $Ap(A)$ is highly indefinite and is much harder for GMRES than the original problem. After $2550$ iterations of GMRES(50), the relative residual almost stalls out at $0.856$.  The vector $b$ appears to have very small components in the eigenvector directions of $A$ that correspond to large eigenvalues.  Thus, the GMRES minimum residual polynomial effectively ignores those large eigenvalues.  

\begin{figure}
\centering
\includegraphics[scale=.6]{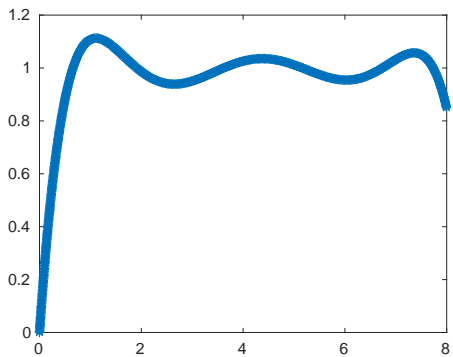}
\caption{Polynomial $\alpha p(\alpha)$ of degree $5$ for the Laplacian generated with $v_0$ as a random vector.}
\label{fig:GoodPolyLapl}
\end{figure}

We now generate a random vector $v_0$ with uniformly distributed elements in $[-1,1]$.  The new polynomial $\alpha p(\alpha)$ of degree $5$ is shown in Figure \ref{fig:GoodPolyLapl}.  The preconditioner works very well; GMRES(50) reaches a relative residual norm of $\num{1.e-8}$ in only $148$ iterations.  The plot of the polynomial shows that the small eigenvalues of $A$ are well-separated and the rest of the spectrum is mapped between $0.8$ and $1.2$.  
It appears that a random vector helps the polynomial to address all parts of the spectrum better than a structured right-hand side vector.  For the remaining experiments in this paper, we let $v_0$ be a random vector.  

\section{Parallel Numerical Results}
\label{sec:parallelresults}
Experiments in this section were performed using the Kodiak cluster at Baylor University.  The cluster has $64$ Cray regular compute nodes, each with dual 18-core Intel E5-2695 V4 (Broadwell) processors and 256GB RAM.  All tests used only one compute node.  

Both examples that follow test finite element discretizations of the convection-diffusion equation $$ -\epsilon \nabla^2 u + \vec{w}\cdot \nabla u = f.$$  The matrices and right-hand sides are generated with Firedrake \cite{Firedrake} software using a function space of continuous piecewise-linear polynomials.
The domain is a 2D unit square mesh centered at the origin with $N=1024$, yielding a matrix of size $n=1050625$.  Similar to Example $6.1.4$ in \cite{ElmanFiniteElts}, $f\equiv 0$ and $$\vec{w} = (2y(1-x^2),-2x(1-y^2)).$$  We use Dirichlet boundary conditions: $u=1$ on boundary $x=1$, and $u=0$ on the remaining boundaries.

We employ GMRES(50), requesting a relative residual tolerance of \num{1e-8} and using two steps of classical Gram-Schmidt orthogonalization.  We generate a random vector $v_0$ and hold it constant for generating all polynomials, regardless of degree or MPI processes.  

\begin{example}
\normalfont
For this example, $\epsilon = 1/2$.  Without preconditioning, GMRES(50) does eventually converge.  The bar graphs in Figure \ref{fig:CoresLowDeg} show solve times over increasing numbers of MPI processes for no preconditioning and polynomial preconditioners of degrees $4$ and $9$.  On one processor, the autodegree selection algorithm chooses degree $9$ as optimal.  The bars are split to show three different timings: 
Time spent in the orthogonalization kernel is indicated by the bottom and middle parts of the bar, for dot products (including norms) and vector updates, respectively.  
The top part of the bar indicates time spent applying $Ap(A)$ to a vector.  Timings for other operations, including polynomial construction, were negligible.  

\begin{figure*}[]
\centering
\subfigure{\includegraphics[scale=.65]{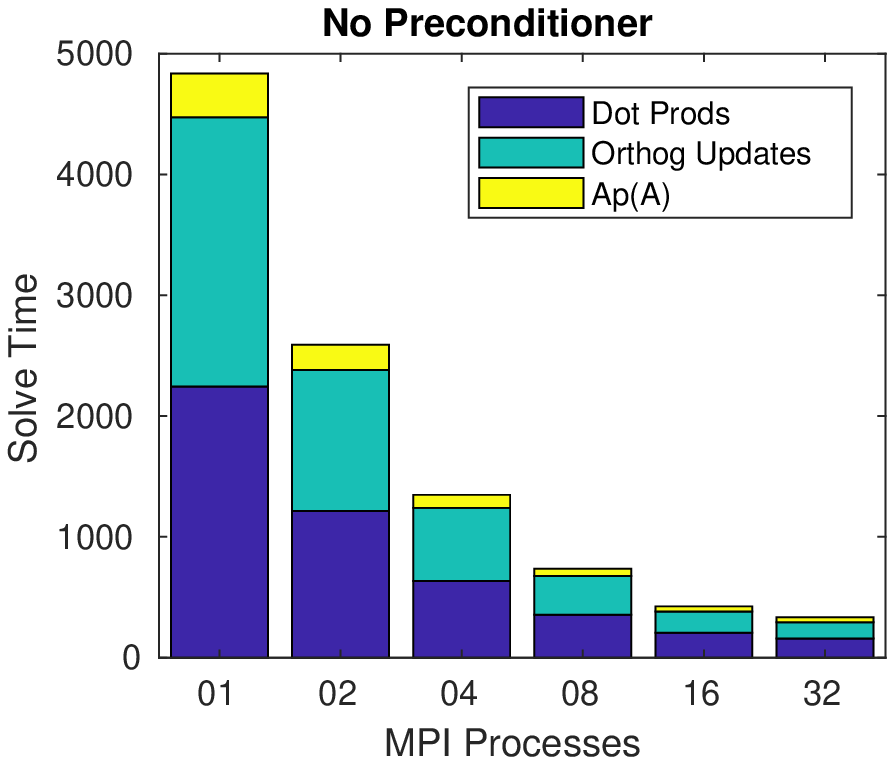}}
\subfigure{\includegraphics[scale=.65]{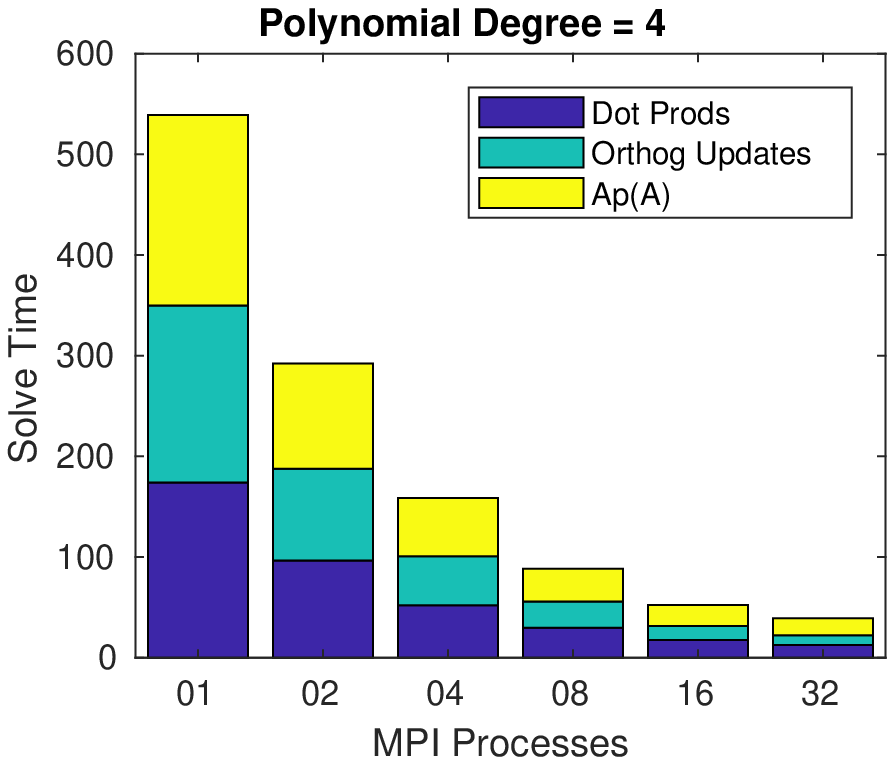}}
\subfigure{\includegraphics[scale=.65]{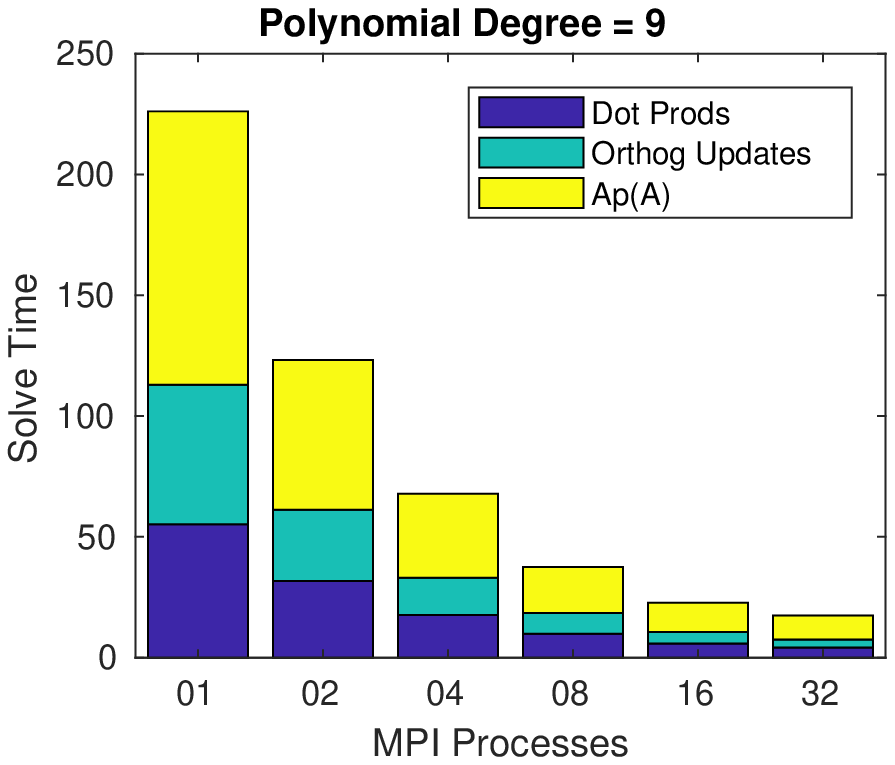}}
\caption{Solve times for convection-diffusion problem using no preconditioning and polynomials of degrees $4$ and $9$ over increasing numbers of MPI Processes. Bottom section of bars give time spent in orthogonalization kernel while top section gives time applying $Ap(A)$.}
\label{fig:CoresLowDeg}
\end{figure*}


Notice first the differences in scaling on the $y$-axes.  The polynomial preconditioner of degree $4$ gives almost $10$ times improvement in solve time over no preconditioning, and degree $9$ gives over $20$ times improvement in solve time over no preconditioning.  Observe also that strong scaling is roughly the same with the preconditioned and unpreconditioned problems. Solve time decreases by about half as we go from $1$ to $2$ MPI processes and by a little less than a half as we add more processes.  

Although running on a single compute node means that communication consists only of reading shared memory, orthogonalization dominates solve time when no preconditioning is used.  In particular, dot products and norms require almost half of the total solve time.  
With a degree $9$ preconditioner, less than one-fourth of the compute time is used for dot products and norms, while a much greater proportion of time is used applying the preconditioned matrix with SpMVs and vector updates.  This shows potential to further reduce solve time by combining polynomial application with a communication-avoiding algorithm such as the Matrix Powers Kernel.  See Section \ref{sec:CA-Alg} for further discussion.  

Surprisingly, there were differences in the polynomial coefficients generated (for fixed degree) with increasing numbers of processors.  Thus, the number of iterations required for convergence varied with the number of MPI processes.  While it would be ideal to have consistent convergence behavior when increasing the number of MPI processes, all the polynomial preconditioners generated here greatly improve convergence.  
\end{example}

\begin{example}
\normalfont 
We modify the convection-diffusion problem from the previous example by choosing $\epsilon = 1/200$.  The increased contribution from the convection term makes this problem too difficult for GMRES(50) to converge without a preconditioner, and polynomial preconditioning alone is ineffective. 
Thus we combine polynomial preconditioning with an ILU preconditioner $M^{-1}$.  We use ILU(0) with no overlap between processors, as implemented in the Trilinos package IFPACK.  We apply ILU preconditioning on the left and polynomial preconditioning on the right. Thus, we are solving the system $M^{-1}Ap(M^{-1}A) y = M^{-1}b$ where $x = p(M^{-1}A)y.$  

Since the ILU preconditioner $M^{-1}$ changes as the number of MPI processes increases (ILU factorizations are computed locally on diagonal blocks of the matrix $A$), the polynomial preconditioners also vary with the number of MPI processes. 
Figure \ref{fig:ILUTimeCores} shows convergence times for polynomial preconditioning of various degrees combined with ILU.  Degree $0$ indicates ILU preconditioning only.  Bar graphs for computations with $1, 8,$ and $32$ MPI processes are shown.  As in the previous example, the bottom and middle sections of each bar indicate orthogonalization time spent on dot products and updates, respectively.  The top section of each bar indicates time spent on all remaining operations, including polynomial and ILU preconditioning. 

\begin{figure*}[]
\centering
\subfigure{\includegraphics[scale=.625]{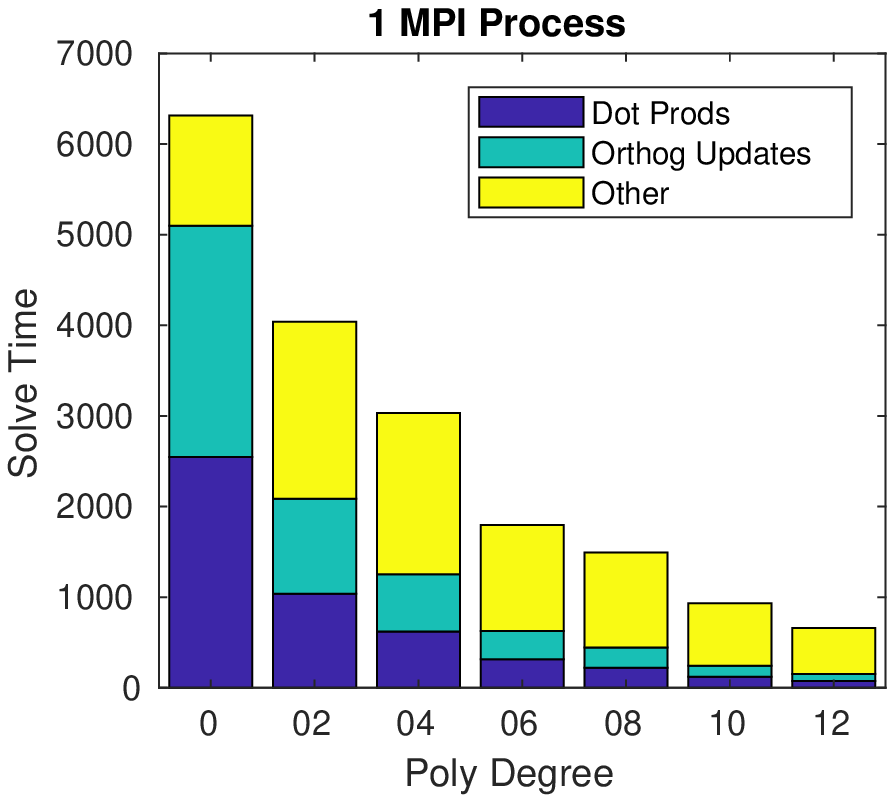}}
\subfigure{\includegraphics[scale=.65]{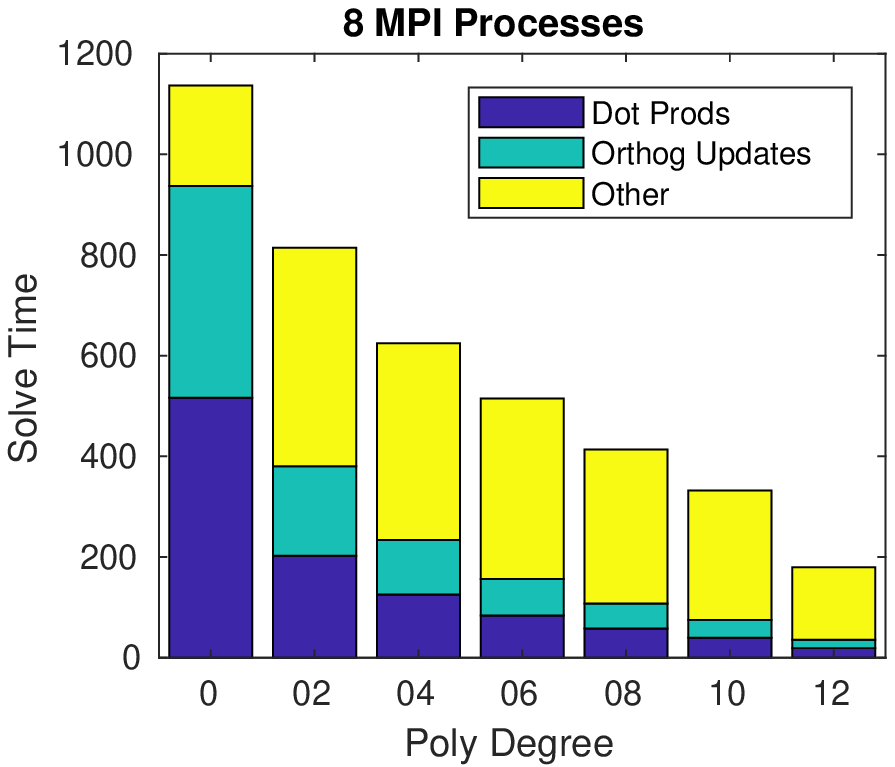}}
\subfigure{\includegraphics[scale=.65]{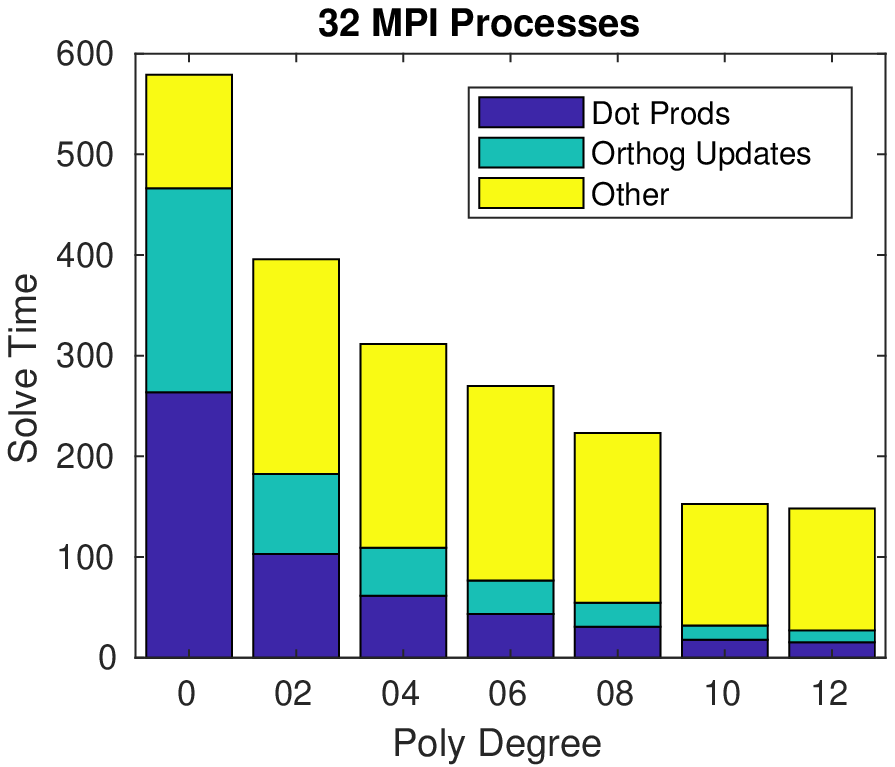}}
\caption{Solve time for polynomial preconditioning combined with ILU over a fixed number of MPI processes.  Degree $0$ indicates ILU preconditioning only.}
\label{fig:ILUTimeCores}
\end{figure*}


Some observations: First, polynomial preconditioning with ILU is significantly better than ILU alone.  Over $1$ MPI process, we attain speedup of almost $10$ times.  Second, this improvement is consistent with increased parallelism.  Even over $32$ MPI processes, we still have almost $4$ times speedup: ILU by itself converges in $576$ seconds, while the degree $12$ polynomial preconditioned GMRES converges in $148$ seconds.  Third, the proportion of time spent in
dot products (orthogonalization) is greatly reduced with polynomial preconditioning.  On $32$ cores with ILU only, dot products and norms consumed $264$ seconds, or about $45\%$ of compute time.  With degree $12$ polynomial preconditioning, they consume only $15$ seconds, or about $10\%$ of compute time.  This suggests that polynomial preconditioning can be a worthwhile addition to ILU and other existing preconditioners. 

%
%
\end{example}

\section{Relation to Communication-Avoiding Methods}
\label{sec:CA-Alg}
Communication-Avoiding Krylov methods, such as CA-GMRES \cite{JAL:HoemmenThesis}, 
are variations of $s$-step Krylov methods. They reduce the number of
communication steps at the cost of more memory and flops. The savings
in global communication and synchronization is important on extreme-scale
parallel systems. There are actually two savings in communication:
(a) Global communication (inner products, orthogonalization)
happens only once every $s$ steps; and (b) with the Matrix Powers
Kernel (MPK), even local communication can be reduced at a cost in memory.
The Matrix Powers Kernel is designed to perform several matrix-vector products consecutively while minimizing reads from memory.  In Communication-Avoiding (CA)-GMRES, the Matrix Powers Kernel is used to form several vectors of a Krylov subspace without orthogonalizing in between.  After all of the SpMVs are performed, then the basis vectors are orthogonalized using the Tall-Skinny QR (TSQR) algorithm.  Unfortunately, CA-GMRES is prone to numerical instability.  It is more stable to form a new Krylov vector from a basis vector that has already been orthogonalized.  The more matrix-vector products are computed before orthogonalization, the more likely the Krylov vectors will begin to lose linear independence.  
Polynomial preconditioned GMRES may provide an avenue for taking advantage of communication-avoiding SpMVs with the Matrix Powers Kernel while avoiding the numerical pitfalls of delayed orthogonalization.  Polynomial preconditioning
can either be used with standard GMRES or within CA methods:
\begin{enumerate}
\item \textbf{Polynomial preconditioned standard GMRES.} We could use the MPK to evaluate the polynomial. Communication occurs as usual in each iteration of the standard GMRES algorithm.
\item \textbf{Polynomial preconditioning within CA-GMRES.} 
Polynomials are ``communication-avoiding'' in the sense that the 
dependency pattern is sparse and it is simple to determine the required 
replication/ghosting of data~\cite{JAL:HoemmenThesis}.
\end{enumerate}
Note that in the first case, we can choose how many powers $t$ we use in 
the MPK, that is, how long to wait between each communication.
The simplest choice is to let $t$ be the degree $deg$ of the polynomial
in (\ref{eq:poly}).
However, this becomes sub-optimal (even impractical) for high degree
polynomials due to the high memory cost. Thus, we are free
to choose $t<deg$.
This is analogous to the fact that in $s$-step methods, 
the length of the MPK, $\bar{s}$, could be different (smaller) 
than $s$.
Observe that polynomial preconditioned GMRES and CA-GMRES will have
essentially the same communication requirements when $t = \bar{s}$
and $deg = s$.
Still, we emphasize they are not equivalent methods.
Our results show that convergence is improved and the number of inner products 
is reduced using polynomial preconditioning.

Also note, by combining polynomial preconditioning and CA-GMRES,
orthogonalization is only needed once every $deg*s$ SpMVs.
Future work includes a more detailed analysis and comparison of
polynomial preconditioned GMRES and CA-GMRES.

\section{Conclusion}
We have shown that polynomial preconditioning can be effective in improving the convergence of GMRES.  Our experiments demonstrate reduction in dot products that helps avoid global communication.  
We showed parallel results on a moderate size cluster. Future work include
experiments on larger problems on highly parallel supercomputers, where
communication is more expensive.  It may also be worthwhile to investigate more stable implementations for constructing the polynomial preconditioner.  

We believe polynomial preconditioning is under-appreciated and
is a good alternative (or complement) to recent communication-avoiding
methods such as CA-GMRES.
It should be made available in high-performance software libraries to
help enable exascale computing.

\section{Acknowledgements}  The first author would like to thank Ron Morgan for several useful discussions and suggestions.  We also thank Rob Kirby for help in generating test problems.  

\bibliographystyle{IEEEtran}
\bibliography{IEEEabrv,JenniferLoe}

\end{document}